\documentclass[12pt]{article}
\usepackage{amssymb}
\usepackage{amscd}
\usepackage[all]{xy}
\usepackage{amsfonts}
\usepackage[centertags]{amsmath}

\newcommand{\be}[1]{\begin{#1}}
\newcommand{\e}[1]{\end{#1}}

\newcommand{\nnd}{\noindent}

\newcommand{\del}{\ensuremath{\delta}}
\newcommand{\T}{\ensuremath{\theta}}

\newcommand{\rh}{\ensuremath{\rho}}
\newcommand{\ga}{\ensuremath{\gamma}}

\newcommand{\al}{\ensuremath{\alpha}}

\newcommand{\csi}{\ensuremath{\xi}}

\newcommand{\sig}{\ensuremath{\sigma}}
\newcommand{\Sig}{\ensuremath{\Sigma}}

\newcommand{\R}{\ensuremath{\mathbf R}}

\newcommand{\Z}{\ensuremath{\mathbf Z}}

\newcommand{\Pro}{\ensuremath{\mathbf P}}
\newcommand{\zd}{\ensuremath{\Z/2\Z}}
\newcommand{\zq}{\ensuremath{\Z/4\Z}}
\newcommand{\zo}{\ensuremath{\Z/8\Z}}
\newcommand{\su}{\ensuremath{S^1}}
\newcommand{\sd}{\ensuremath{S^2}}

\newcommand{\spin}{\ensuremath{Spin}}

\newcommand{\ie}{i.\,e.}

\newcommand{\tM}{\ensuremath{\tilde{M}}}

\newcommand{\tK}{\ensuremath{\tilde{K}}}
\newcommand{\tN}{\ensuremath{\tilde{N}}}
\newcommand{\tH}{\ensuremath{\tilde{H}}}

\newcommand{\tSig}{\ensuremath{\tilde{\Sigma}}}
\newcommand{\tf}{\ensuremath{\tilde{f}}}

\newcommand{\ttau}{\ensuremath{\tilde{\tau}}}

\newcommand{\tc}{\ensuremath{\tilde{c}}}

\newcommand{\cl}[1]{\ensuremath{[{#1}]}}

\newcommand{\Imm}{\ensuremath{Imm_{\xi}(F,M)}}

\title{Cobordism of immersions of surfaces in non-orientable 3-manifolds.}

\author{Rosa Gini}

\date{January 2000}

\begin{document}

\maketitle

\begin{abstract}Some properties of non-orientable 3-manifolds
are shown. The semi--group of cobordism of immersions of surfaces
in such manifolds is computed and proven actually to be a group.
Explicit invariants are provided.

\end{abstract}
\section*{Introduction}


Following the definitions of \cite{WelCGI}, \cite{PinRHC},
\cite{BeSSPS} we will say that two immersions $f$ and $f'$ defined
on compact closed surfaces not necessarily connected $F$ and $F'$
and taking values in the same 3-manifold $M$ are  {\em cobordant}
if there exists a cobordism $X$ between $F$ and $F'$ and an
immersion $\Phi$ of $X$ in $M\times I$ that restricts to $f\times
\{0\} $ and to $f'\times\{1\}$. Once fixed the manifold $M$ the
set $N_2(M)$ of cobordism classes of immersions of surfaces in $M$
is a semi-group with the composition law given by disjoint union.
That this be a group when $M=\R^3$ is given a priori by the fact
that inverses are provided by composition with a reflection in a
plane. In fact $N_2(\R^3)$ is isomorphic to $\Z/8\Z$, as is proved
in \cite{WelCGI}; explicit invariants are given in \cite{BroGKI}
and \cite{PinRHC}; a generator is the so-called {\em right
immersion of Boy} of the projective space, see \cite{PinRHC}. That
$N_2(M)$ be a group when $M$ is a generic manifold is not
straightforward; it is proven in \cite{BeSSPS} for  orientable
3-manifolds  that $N_2(M)$ is the finite set $H_2(M,\zd)\times
H_1(M,\zd)\times\zo$ with a composition law that twists the
compositions of the factors. We adapt here their constructions to
the non-orientable case, and obtain in theorem~\ref{cobordismo}
that again $N_2(M)$ is a finite group, with support the set
$H_2(M,\zd)\times H_1(M,\zd)\times\zd$ and a composition law
similar to the one of the orientable case. The crucial point that
causes a non-orientable 3-manifold to have a ``smaller" group than
an orientable 3-manifold with same homology is
proposition~\ref{C(M)}, and comes from the fact that in a
non-orientable environment there exist isotopies that reverse
orientation.

The first section is essentially a review of results of
\cite{HaHIST}, \cite{BeSSPS} and \cite{KiTPSL}; some remarks and
properties are original and are used in the following; in
particular in proposition~\ref{teoisotropia} we compute explicitly
the isotropy group for the action of adding kinks introduced in
\cite{HaHIST}, and in theorem~\ref{strisce} we classify bands
(that is, immersions of annuli or Moebius bands) in a
non-orientable 3-manifolds up to regular homotopy. In the second
section we develop the computation of the cobordism group.

This paper is an extended version of the talk given in Palermo in
September 1999, during the Congress {\em ``Propriet\`a Geometriche
delle Variet\`a Reali e Complesse: Nuovi Contributi Italiani"}.

\tableofcontents

\section{Some properties of non-orientable 3-manifolds.}

From now on a loop $c$ in a non-orientable manifold $X$ will be
said to be {\em orientable in $X$ }if it preserves the orientation
of $X$, that is, if $TX_{|c}$ is orientable, and will be said to
be {\em non-orientable in $X$ }otherwise.

\subsection{Projective and anti-equivariant framings.}

It is well-known that orientable 3-manifolds are parallelizable. A
non-orientable 3-manifold obviously is not, but its tangent bundle
is still as simple as a non-orientable vector bundle can be, that
is, its structure group can be reduced to $\{1,-1\}$. A calculus
of characteristic classes in \cite{HaHIST} proves in fact that if
$M$ is a non-orientable 3-manifold then
\[
TM\cong \det M\oplus  \det M\oplus \det M.
\]

This means that it is possible to define on $M$ a {\em projective
framing}, i.e.\ a triple of linearly independent vector fields,
well-defined up to sign. A triple $\{v_1,v_2,v_3\}$ of linearly
independent vector fields on  the orientation covering space \tM\
of $M$ is  said to be an {\em anti-equivariant framing} if
\[
v_i(\sig (P))=-\sig_*(v_i(P)) \hspace{1cm} \forall P\in\tM,
i=1,\dots,3,
\]
where \sig\ is the covering translation of \tM; an
anti-equivariant framing projects on a projective framing. Given a
projective framing of $M$ one can always find an anti-equivariant
framing of \tM\ that projects on it: we call such an
anti-equivariant framing a {\em lifting} of the projective
framing.

We will assume, from now on, that all framings are orthonormal
with respect to a fixed metric.

\subsection{Isotropy in adding kinks.}\label{isotropia}

\begin{figure}
\begin{picture}(270,150)
\put(30,30){ \includegraphics{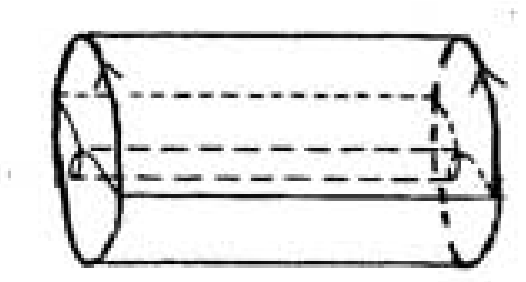} }
\end{picture}
\caption{The local modification of adding a kink.}\label{addkink}
\end{figure}

The main result of \cite{HaHIST} is to set up an action $*$  of
$H^1(F,\zd)$ on the set $Imm_{\csi}(F,M)$ of regular homotopy
classes of immersions in a homotopy class \csi\ of immersions of a
surface $F$ in a 3-manifold $M$. This action is explicitly
described in a beautiful geometric way: if $f$ is a representative
of an element of \Imm\ and \al\ is in $H^1(F,\zd)$ then a
representative of $f*\al$ is obtained by modifying $f$ in a
tubular neighborhood $N$ of a dual curve to \al, in a way the
authors describe as {\em adding a kink}: in case both of the curve
in $F$ and its image in $M$ are orientable this local modification
is the (local) composition of $f$ with the immersion of an annulus
in a solid torus shown in figure~\ref{addkink}; in the other cases
the modification is similar.
If $c$ is a dual curve to \al\ then $f*\al$ will also be denoted
by $f*c$.

The action of adding kinks is proven to be transitive in any
homotopy class. Isotropy depends on a property of the class: we
say a homotopy class \csi\ is {\em odd} if both of $F$ and $M$ are
non orientable and if there exists a self-homotopy $H$ of a map
$f\in\xi$ such that $H(x,-)$ is a loop non-orientable in $M$ for
$x\in F$; we call {\em odd} such a homotopy; the action of
$H^1(F,\zd)$ on $Imm_{\csi}(F,M)$ has then a group of isotropy of
order 2 in any point if \csi\ is odd. If \csi\ is not odd we say
it is {\em even}, and isotropy is trivial in this case. This
means, in particular,  that the set $Imm_{\csi}(F,M)$ has the
cardinality of $H^1(F,\zd)$ if \csi\ is even, and half its
cardinality if \csi\ is odd.

We compute explicitly the isotropy group for the action on odd
classes.

\be{prop}\label{teoisotropia}Let \csi\ be a odd homotopy class of
immersions of a non-orientable surface $F$ in a non-orientable
3-manifold $M$. Then the isotropy group of any regular homotopy
class in \csi\ is the vector subspace of $H^1(F,\zd)$ generated by
$w_1(F)$. \e{prop}

\be{dimostraz} The action of $H^1(F,\zd)$ on \Imm\ is not
explicitly defined in \cite{HaHIST}. What is stated there is that,
given a point $f$ of \Imm, it is possible to define a
correspondence $C_f$ between \Imm\ and $H^1(F,\zd)$ which is 1-1
if \csi\ is even and 1-2 if \csi\ is odd; remark that $C_f^{-1}$
is then always a function. This action is explicitly defined in
\cite{BeSSPS} as
\[
\be{array}{cccl}*:& \Imm \times H^1(F,\zd)&\longrightarrow&\Imm\\
&(f,\al)&\mapsto& \al*f:=C_f^{-1}(\al).\e{array}
\]

We prove that, if \csi\ is odd, then $C_f(f)=\{1,w_1(F)\}$. The
correspondence $C_f$ is constructed in several steps. The first is
to associate to every $g\in\Imm$ the homotopy class of its
differential: that this be a bijection is a consequence of
\cite{HirIM}; to $f$ is then associated $df$. The second step is
where non-injectivity appears: to each differential one associates
two elements of the set $Bun_f(TF,TM)$ of homotopy classes of
bundle map that commute with $f$, via the choice of a projective
framing of $M$. It follows from the construction given in
\cite{HaHIST} that the classes associated to $df$ are $df$ itself
and $-df$. The following steps give a bijection between
$Bun_f(TF,TM)$ and $H^1(F,\zd)$, that associates the identity to
$df$; we call $w$ the element of $H^1(F,\zd)$ associated to $-df$:
we are then left to prove that $w=w_1(F)$.

Let $c$ be a curve on $F$. We describe how to determine $w(c)$.
 Consider on
$c$ the fiber bundle $\tau:=f^*(TM)=TF_{|c}\oplus\nu_c$ and the
automorphism of $\tau$ given by $-1_{TF_{|c}}\oplus 1_c$. If
$f(c)$ is orientable in $M$ then the projective framing determines
two opposite framings of the trivial bundle $\tau$; the map
$-1_{TF_{|c}}\oplus 1_c$ read in one of these framings gives a
closed path in $SO(3)$; $w(c)$ is then defined to be the class of
this path in $H_1(SO(3),\zd)=\zd$, that does not depend on the
choice between the two framings (nor on the choice of $c$ between
representatives of its class modulo 2) . If $f(c)$ is not
orientable then $\tau$ is not trivial; the projective framing of
$M$ restricted to $\tau$ defines two triples of orthonormal vector
fields discontinuous in a point; choose one of them
$\{v_1,v_2,v_3\}$; define a path of $SO(3)$ by
\[ P\mapsto\text{\footnotesize linear transformation between
$\{v_1(P),v_2(P),v_3(P)\}$ and $-1_{TF_{|c}}\oplus
1_c(v_1(P),v_2(P),v_3(P))$,}
\]
for all $P\in c$; thought the vector fields are discontinuous this
path is well-defined, continuous and closed, and its homology
class is independent of the choice between the two opposite triple
of discontinuous vector fields. The value of $w(c)$ is this
homology class, again considered as an element of \zd.

We must now compute these values. Recall that a point of $SO(3)$
is identified by a point of \sd, which indicates the oriented
rotation axis, and a number between 0 and $\pi$, which indicates
the rotation angle taken in the positive sense given by
orientation; the facts that the null rotation around any axis is
the identity and that two rotations of $\pi$ along opposite axes
coincide give a bijection between $SO(3)$ and the projective space
$\Pro^3$, which is in fact a well-known homeomorphism. The point
in $SO(3)$ associated to a point $P$ of $c$ is the rotation of
$\pi$ whit axis any of the two normal vectors in $f(P)$ to $f(N)$,
$N$ being a tubular neighborhood of $c$ in $F$. We compute the
homology of paths of this sort as the intersection number with a
generator of $H_1(SO(3),\zd)$: we take the path of rotations of
$\pi$ with rotation axis rotating itself of $\pi$ in the $xz$
plane, from $(1,0,0)$ to $(-1,0,0)$.

Consider first the case that $f(c)$ is orientable in $M$, that is
$\tau$ is framed by the projective framing. We then consider a
homeomorphism of the standard torus in $\R^3$ to a tubular
neighborhood of $f(c)$ which sends the standard framing to the
framing of $\tau$. The image of $f(N)$  results as a band in this
standard torus, and the normal vector to this bands is parallel to
the $xz$ plane (in points where this intersection may be
considered transverse) a number of times which is even if the band
is orientable and odd if the band is a Moebius band: so we have
shown that $w(c)$ coincides with $w_1(F)(c)$ when $f(c)$ is
orientable in $M$.

We are then left to the case $f(c)$ non-orientable in $M$. We may
suppose to have fixed a lifting to \tM\ of the projective framing
of $M$. Consider the non trivial double covering \tN\ of $N$, let
\tf\ be a map of this covering to \tM\ that commutes with $f$
restricted to $N$, let \ttau\ be $\tf^*(TM)$: this is a trivial
bundle framed by the fixed framing of \tM. Again we fix a
homeomorphism of the standard torus with standard framing to a
tubular neighborhood of $\tf(\tc)$ in \tM\ with the fixed framing;
the image of the band $\tf(\tN)$ in this torus gives count two
times of the band $f(N)$, so consider its (transverse)
intersection with, say, the semi-space $y\geq0$; the number we are
looking for is the number of times, modulo 2, that the normal
vector to this part of the band is parallel to the $xz$ plane
(again this intersection can be taken to be transverse); again
this number is even if $N$ is orientable and odd if $N$ is a
Moebius band, so that again $w(c)=w_1(F)(c)$, and this ends the
proof.~$\triangle$

\e{dimostraz}


\subsection{Bands in orientable 3-manifolds.}

The meaning of the main theorem in \cite{HaHIST} is that the
regular homotopy class of a map in its homotopy class is
determined by the behaviour of the map in tubular neighborhoods of
the curves of the surface. These tubular neighborhoods are either
annuli or Moebius bands; we call {\em band in a 3-manifold} the
immersion of an annulus or a Moebius band. In \cite{BeSSPS} the
r\^ole of bands in this subject becomes more evident, and their
properties are also a necessary brick in the computation of
cobordism in the case of orientable 3-manifolds.

We review the properties of bands in orientable 3-manifolds, then
we look at the case of bands in a non-orientable 3-manifold.
\subsubsection{Bands in $\R^3$.}

Consider an embedded band in $\R^3$, and consider the linking
number of its core with its (possibly non-connected) boundary. We
think of this linking number in the following way: we pick on the
boundary of a tubular neighborhood $N$ of the core of the band a
preferred basis for $H_1(\partial N,\Z)$, with the orientation of
the longitude coherent with an orientation of the core; we then
give to the meridian the orientation that makes the global
orientation off $\partial N$ compatible with the orientation of
the environment; the intersection of the band with $\partial N$
represents a homology class, and the linking number is the
meridian coordinate of this class in the preferred basis. This
number is even if the band is orientable and odd if it is a
Moebius band, so that its class modulo 2 is a total topological
invariant of the band.

\begin{figure}
\begin{picture}(270,150)
\put(30,30){ \includegraphics{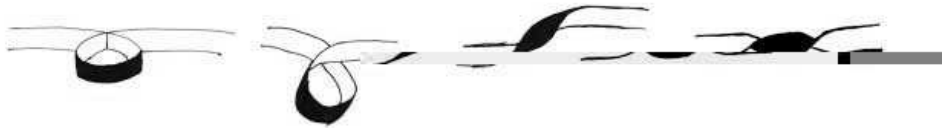} }
\end{picture}
\caption{Two half twists become a kink.}\label{ricciolo}
\end{figure}

\begin{figure}
\begin{picture}(270,150)
\put(30,30){ \includegraphics{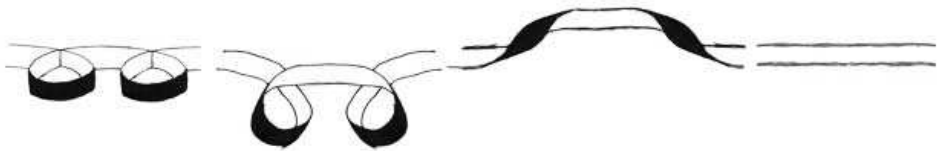} }
\end{picture}
\caption{Two kinks become no kinks.}\label{riccioli}
\end{figure}

Given any band in $\R^3$ there exists a regular homotopy between
its core and the standard circle in the $xy$ plane in $\R^3$;
extend it to a regular homotopy between a tubular neighborhood of
the core of the band and the standard torus of $\R^3$: we obtain a
band regularly homotopic to the original and that differs from the
annulus in the $xy$ plane by a number of half twists given by the
linking number. Now if we isolate 4 half twists there exists a
regular homotopy  relative to the rest of the band which
transforms 4 half twists in a couple of kinks and then in a piece
of band with no twisting, see figures~\ref{ricciolo}
and~\ref{riccioli};
this local construction reduces by regular homotopy the original
band to one of four models, classified by the linking number
modulo 4. This number is in fact an invariant in \zq\ of the
regular homotopy class of the band (see~\cite{GuMETR}, page 114),
so it coincides with the original linking number (modulo 4). This
proves the following:

\be{prop} In $\R^3$ the linking number modulo 4 between the core
of a band and the boundary oriented in a coherent way is a total
invariant of the  regular homotopy class of the band.~$\triangle$
\e{prop}

\vspace{0.5cm}

We call this invariant in \zq\ {\em number of half twists of the
band}.

\subsubsection{\spin\ structures and preferred longitudes in orientable
3-manifolds.}\label{framingdelcerchio}

The problem in extending the definition of half twists of a band
to generic 3-manifolds is that in a 3-manifold there is no
preferred basis for the homology of the boundary of the tubular
neighborhood of a knot: a meridian is still defined, and the
orientation of the environment can give it an orientation as it
does in $\R^3$, but there is no way, a priori, to distinguish
between longitudes. Remark that the local modification of
figure~\ref{riccioli}  can be performed in any manifold, so that
up to regular homotopy we are still interested in the class of
rest modulo 4 of the integer coordinate; this means that it is
enough to choose a homology class of longitudes modulo 2, since
longitudes that differ by two meridians give integer coordinates
belonging to the same class modulo 4 (remark that the boundary of
a band covers the core twice). In \cite{KiTPSL} it is shown how to
use a \spin\ structure on an orientable 3-manifold to determine
this choice, that reduces to the classical one when $M=\R^3$ with
its unique \spin\ structure. We give here a slightly different way
to  define the same choice.

Orientable vector bundles on a circle are always trivial. Framings
of an oriented vector bundle of rank 3 on a circle are 2, up to
homotopy. This is because a homotopy class of framings can be
considered as a homotopy class of sections of the associated
principal bundle, which, being trivial, is homeomorphic to
$\su\times SO(3)$: and $\pi_1(\su\times SO(3))$ contains four
classes, only two of which can be realized as sections of the
bundle.

Consider \su\ embedded in $\R^3$ in the standard way, and consider
on it the fiber bundle $\tau$ given by restriction of the tangent
bundle to $\R^3$. Define on $\tau$ two framings containing in each
point the tangent vector to \su: let $e_0$ be the oriented framing
containing the vector parallel to the $z$-axis as second element;
let $e_1$ be the framing containing in the point
$(\cos\T,\sin\T,0)$ of \su\ the vector
$\cos\T(0,0,1)+\sin\T(\cos\T,\sin\T,0)$ as second element. These
two framings are not equivalent, since the curves they describe in
the principal bundle of $\tau$ differ by a generator of the
homotopy group of $SO(3)$: this means that any other framing of
$\tau$ is equivalent either to $e_0$ or to $e_1$. In particular
the standard framing is homotopic to $e_1$, since $e_1$ extends to
a framing of the tangent bundle to $\R^3$ restricted to the disc,
by $$ (\rho\cos\T,\rho\sin\T,0)\mapsto((-\sin\rho\T,\cos\rho\T,0),
(\sin\rho\T\cos\rho\T,\sin^2\rho\T,\cos\rho\T ),\text{third
orthonormal}), $$ and remark that the framing of the trivial
bundle on the disc is unique up to homotopy, the disc being
contractile. But now the second vector of either $e_0$ or $e_1$
frames the normal bundle to \su\ in $\R^3$, and in particular
detects a longitude;
so we get close to what we are looking for. Remark that the
longitude which is preferred in the classical definition (that is,
the one having linking number 0 with the core) is the one given by
$e_0$, and that the longitude detected by $e_1$ belongs to the
other class.

On the other side, \spin\ structures on an oriented vector bundle
on a manifold $B$ are acted on simply transitively by
$H^1(B,\zd)$, so that \spin\ structures on an oriented vector
bundle of rank 3 on the circle are again 2. It is also true that a
framing of a trivial bundle induces naturally a \spin\ structure:
in fact a \spin\ structure on a vector bundle can be defined as a
double covering of the associated principal $SO(n)$-bundle which
be non trivial when restricted to the fiber (see~\cite{MilSSM}),
and a framing, giving a bundle equivalence between the principal
bundle and $SO(n)\times B$, induces such a covering by pull-back
of the standard covering by $\spin(n)\times B$:
\[
\xymatrix{\text{induced spin structure}\ar@{.>}[r]\ar@{.>}[d]
&\spin(n)\times B\ar[d]\\ \text{principal bundle
}\ar[r]_{\text{framing}}&SO(n)\times B. }
\]
The definition of \spin\ structure via double covering implies
that \spin\ structures can be described by the elements of
$H^1(\text{principal bundle},\zd)$ that restrict to the fiber as
generators of $H^1(SO(n),\zd)$ . Two \spin\ structures are {\em
equivalent} if the two double coverings are equivalent, \ie\ if
the associated cohomology classes are the same.

Now, $\pi_1(\su\times SO(3))$ and $H^1(\su\times SO(3),\zd)$ are
isomorphic, and the correspondence we gave between homotopy
classes of sections of the principal bundle and cohomology classes
that restrict to a fiber as a generator is bijective in this case,
so that speaking of \spin\ structures or of framings of an
oriented vector bundle of rank 3 on \su\ is the same thing.

We are now ready for the definition. Let $M$ be an oriented
3-manifold with a fixed \spin\ structure, let $K$ be a knot  in
$M$; the \spin\ structure restricted on $TM_{|K}$ gives a framing;
choose a diffeomorphism of the standard \su\ in $\R^3$ with
standard framing to $K$ in $M$ with this framing, deform by a
homotopy the standard framing to $e_1$ and pull this deformation
back to $M$; the second vector of the framing, with usual
identification of the normal bundle with the tubular neighborhood
$N$ of the knot, gives now a curve on $\partial N$; the {\em
preferred longitude} on $\partial N$ will be the homology class in
$H_1(\partial N,\zd)$ that doesn't contain this curve. It is
straightforward to verify that this choice doesn't depend on the
diffeomorphism nor on the deformation: in fact an oriented framing
of $T\R^3$ restricted to the standard \su\ and containing the
tangent vector to \su\ is homotopic either to $e_0$ or to $e_1$
according to the homology class modulo 2 of the curve they pick on
the standard torus is equal to the class to the standard longitude
or not. This proves the following characterization:

\be{prop} Any deformation of the framing of $TM_{|K}$ that takes
it to a framing containing the tangent vector to $K$ chooses a
curve in $\partial N$ which is the non-preferred
longitude.~$\triangle$ \e{prop}

 \vspace{0.5cm}

This description of the preferred longitude is particularly useful
when the \spin\ structure of the manifold is itself associated to
a framing of the whole manifold. In this case to assign the
preferred longitude it is even not necessary to mention the \spin\
structure.


We prove that our definition coincide with the definition of {\em
even framing} in \cite{KiTPSL}, page~209:

\be{lem}Let $M$ be an oriented 3-manifold with a fixed \spin\
structure, and let $K$ be a knot in $M$; then any even framing of
the normal bundle to $K$ represents a preferred longitude.

 \e{lem}

\be{dimostraz} The choice of even framing only depends on the
\spin\ structure of $TM$ restricted to $K$, so  it coincides
either with the choice of preferred longitude or with its
opposite. But the two choices coincide in the case of the standard
circle in $\R^3$, where they both assign the class of the standard
longitude of the torus, so they coincide in every
case.~$\triangle$

\e{dimostraz}

If the knot represents the trivial class in $H_1(M,\zd)$  the
preferred longitude does not depend on the choice of the \spin\
structure, since any two structures differ by the action of an
element of $H^1(M,\zd)$, that  doesn't affect this knot. In this
case the preferred longitude can be realized in a geometric way:

\be{prop}If a knot $K$ represents the trivial class in
$H_1(M,\zd)$ then take any embedded surface $F$ such that
$\partial F=K$; the (transverse) intersection of $F$ and $\partial
N$ is a preferred longitude.
 \e{prop}

\be{dimostraz} This property is the content of theorem~4.3 in
\cite{KiTPSL}.~$\triangle$

\e{dimostraz}

\be{coroll} If $M=\R^3$ the choice of preferred longitude reduces
to the classical.~$\triangle$ \e{coroll}

\vspace{0.5cm}

Remark that the definition of preferred longitude is invariant
under regular homotopy. For more details on this subject see
\cite{KiTPSL}.

In \cite{BeSSPS} the definition of even longitude allows to extend
the notion  of number of half twists  to bands in generic oriented
3-manifolds. What is (implicitly) obtained is that given a knot
$K$ there are 4 regular homotopy classes of bands having $K$ (or a
knot homotopic to it) as core, and that  the number of half twists
is a total invariant.

\subsection{Bands in non-orientable 3-manifolds.}

\subsubsection{Equivariant \spin\ structure.}

The first  and less expensive attempt to extend the definition of
half twists to bands immersed in a non-orientable 3-manifold $M$
is to look only at bands whose core is orientable in $M$, that is,
bands that have two homeomorphic preimages in the orientation
double covering \tM\ of $M$: the {\em number of half twist} of the
band can be defined as the number of half twists of one of its
preimages, provided we fix on \tM\ a \spin\ structure that makes
the choice between the two preimages inessential. Call \sig\ the
covering translation of \tM: we say a \spin\ structure on \tM\ is
{\em equivariant} if, given an embedded circle $K$ in \tM\ such
that $\sig K\cap K=\emptyset$, whenever $l$ is a curve on the
boundary of a tubular neighborhood of $K$ that represents a
preferred longitude, then $\sig l$ represents a preferred
longitude.

\be{lem} The \spin\ structure on \tM\ induced by an
anti-equivariant framing is equivariant.

\e{lem}

\be{dimostraz} Remark that \sig\ restricted to a tubular
neighborhood $N$ of $K$ does not preserve the homotopy class of
the framing, since \sig\ is orientation reversing. Call
$\{v_1,v_2,v_3\}$ the given anti-equivariant framing: the
longitude $l$ belongs to the class which doesn't contain the curve
detected on $\partial N$ by a deformation of
$\{v_1,v_2,v_3\}_{|K}$; the same deformation composed with \sig\
gives a deformation of $\{\sig v_1,\sig v_2,\sig v_3\}_{|\sig
K}=\{-v_1,-v_2,-v_3\}_{|\sig K}$, which picks on the boundary of
the tubular neighborhood of $\sig K$ the longitude $\sig l$.

We are then reduced to show that the curve on the standard torus
in $\R^3$ detected by the non-oriented framing $-e_1$ belongs to
the class which doesn't contain the standard longitude (recall we
are talking of classes modulo 2). But this is evident.~$\triangle$

\e{dimostraz}

Now put on \tM\ an equivariant \spin\ structure. Let \Sig\ be a
band in $M$ whose core is orientable in $M$, let \tSig\ be one of
its liftings to \tM; the other is \sig\tSig. We have

\be{prop}\label{opposto}The number of half twists of \tSig\ is
opposite to the number of half twists of \sig\tSig.

\begin{figure}
\begin{picture}(270,150)
\put(30,30){
\includegraphics{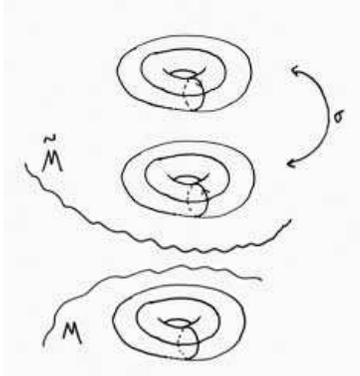} }
\end{picture}
\caption{The covering of the tubular neighborhood of a curve
orientable in $M$.}\label{cuoreorientabile}
\end{figure}

\e{prop}

\be{dimostraz} Call $K$ the core of \Sig, with a fixed
orientation, and call $N$ a tubular neighborhood of $K$, see
figure~\ref{cuoreorientabile}.
 Let
$\{m,l\}$ be a preferred basis for $H_1(\partial N,\Z)$: recall
the orientation of $l$ is chosen according to the orientation of
$K$, and the orientation of $m$ is the one that, together with the
orientation of $l$, gives to $\partial N$ the orientation coming
from \tM. This choice makes the meridian coordinate of a curve in
$\partial N$ independent of the choice of orientation for $K$.

Recall that \sig\ is orientation reversing. This means that
$\{\sig m,\sig l\}$ is not a preferred basis for the homology of a
tubular neighborhood of $\sig K$, whereas $\{-\sig m,\sig l\}$ is
one. This gives the thesis.~$\triangle$

\e{dimostraz}

But an even number and its opposite are congruent modulo 4; this
proves:

\be{coroll}The number of half twists is well defined for
orientable bands whose core is orientable in $M$.~$\triangle$

\e{coroll}

\subsubsection{Regular homotopies between bands in non-orientable
3-manifolds.}

In trying to extend further the definition of number of half
twists to generic bands in non-orientable 3-manifold we look at
regular homotopies in a non-orientable environment: we will see
that in a non-orientable environment there are ``more" regular
homotopies. First recall the notion of odd self-homotopy in a
non-orientable manifold $X$ (see \S~\ref{isotropia}): a
self-homotopy $H$ is odd if the closed path $H(x,-)$ is
non-orientable in $X$ for any $x$.

\be{teo}\label{strisce}Let $M$ be a non-orientable 3-manifold,
and let $K$ be an embedded circle
 in $M$; then

\be{enumerate}

\item if $K$ is orientable in $M$ and doesn't admit any odd self-homotopy
then there are 4 regular homotopy classes of bands having a circle
homotopic to $K$ as core;

\item if $K$ is orientable in $M$ and admits an odd
self-homotopy then there are 3 regular homotopy classes of bands
having a circle homotopic to $K$ as core, in particular any two
non-orientable bands of this type are regularly homotopic;

\item if $K$ is non-orientable in $M$ then there are 3 regular homotopy classes
of bands having a circle homotopic to $K$ as core, in particular
any two non-orientable bands of this type are regularly homotopic;
the two classes of orientable bands differ by a reparametrization.

\e{enumerate}

\begin{figure}
\begin{picture}(270,150)
\put(30,30){ \includegraphics{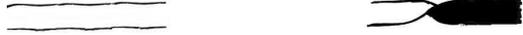} }
\end{picture}
\caption{Adding a local twist.}\label{localtwist}
\end{figure}

\e{teo}

\be{dimostraz}Let \Sig\ be a band whose core is homotopic to $K$.
Take the core of \Sig\ to $K$ with a regular homotopy and extend
it (for example via an exponential map) to a tubular neighborhood:
this takes \Sig\ to a band regularly homotopic and having $K$ as
core. We now think of \Sig\ as a band having $K$ as core. By means
of the local modification shown in figure~\ref{riccioli} we can
reduce any band having $K$ as core to one of four models; call
$\Sig_0$ one of the models of orientable band, the others differ
from it only locally: call {\em adding a local twist} the
operation of substituting to a piece looking as in the left side
of figure~\ref{localtwist}
a piece looking as in the right side of the same figure; $\Sig_i$
has $i$ local half twists, $i$ being 1 or 2, $\Sig_{-1}$ has one
local half twist in the opposite direction to the half twist of
$\Sig_1$. We have to decide under which conditions these models
can or cannot admit regular homotopies. Recall that  bands with
number of half twists that differ modulo 2 are defined on
non-homeomorphic domains, so that the notion of regular homotopy
is only possible between $\Sig_0$ and $\Sig_2$ and between
$\Sig_1$ and $\Sig_{-1}$.

Let now $K$ be a curve orientable in $M$. Its tubular neighborhood
$N$ is then a solid torus. The embedding $K$ lifts to \tM\ in two
embeddings \tK\ and \sig\tK; any homotopy between the two liftings
projects to an odd self-homotopy of $K$, any self-homotopy of one
of the liftings projects to an even self-homotopy of $K$. The
liftings of the embedding of $N$ are two solid tori, and the
covering projection restricted to either of these tori gives
opposite orientations to $N$, see figure~\ref{cuoreorientabile}.
If we extend a homotopy between \tK\ and \sig\tK\ to tubular
neighborhoods it projects to an odd self-homotopy of $N$ that
reverses the orientation, if we extend a self-homotopy of either
\tK\ or \sig\tK\ to a tubular neighborhood it projects to an even
self-homotopy that preserves the orientation. On the other side
any self-homotopy of $K$ lifts either to a self-homotopy of \tK,
if it is even, or to a homotopy between \tK\ and \sig\tK, if it is
odd.

Now if $K$ does not admit any self-homotopy (for example if \tK\
and \sig\tK\ are not homotopic) then any regular homotopy between
two different models would lift to a regular homotopy between
bands having \tK\ as core and different number of half twists, and
this is not possible by the classification of bands in orientable
manifolds: this proves the first part of the theorem.

If $K$ admits an odd self-homotopy then take a regular odd
self-homotopy and extend it to a regular self-homotopy $H$ of $N$.
Call $\tSig_i$ the lifting of $\Sig_i$ to a band having \tK\ as
core; we might suppose that $\tSig_i$ has $i$ half twists. The
lifted regular homotopy \tH\ takes $\tSig_i$ to a band having
\sig\tK\ as core, and the same number of half twists. Recall that,
by \ref{opposto}, $\sig\tSig_i$ has $-i$ half twists, so, if $i$
is odd it is not the final image of \tH. But $\sig\tSig_i$
projects to $\Sig_i$, so that the final image of $\Sig_1$ under
$H$ is the other model, that is, $\Sig_{-1}$. The same
construction applied to $\Sig_0$ shows that a odd homotopy
preserves the model when the band is orientable, and that even
homotopy preserve comes from the first part. This proves the
second assertion.

\begin{figure}
\begin{picture}(270,150)
\put(30,30){ \includegraphics{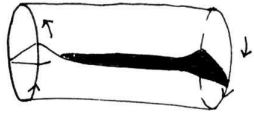} }
\end{picture}
\caption{The twisting of a band with core non-orientable in
$M$.}\label{caramella}
\end{figure}

We are left to the case when $K$ is a non-orientable curve in $M$.
Its tubular neighborhood $N$ is then a solid Klein bottle. We can
fix a diffeomorphism of $N$ to the following model of solid Klein
bottle: consider $D^2\times I$ modulo the relation
$(P,0)\sim(\rho(P),1)$, \rh\ being the reflection of $\R^2$ in the
$x$ axis restricted to $D^2$. We can consider that $\Sig_0$ is
$[-1,1]\times\{0\}\times I/\sim$. Remark now that, when we rotate
$[-1,1]\times\{0\}\times\{0\}$ in the positive direction, the
other copy of this segment, that is
$[-1,1]\times\{0\}\times\{1\}$, must rotate in the negative
direction, see figure~\ref{caramella}.
The effect on the band is that, after a rotation of $\pi$,
$\Sig_0$ becomes a band having the same support as $\Sig_2$, and
differing from it only by a reparametrization. Remark that the two
boundary components of $\Sig_0$ belong to different homology
classes.

On the other side if we apply this isotopy to $\Sig_1$ and rotate
of $\pi/2$ we obtain a band with support $\{0\}\times[-1,1]\times
I/\sim$; if we apply the opposite rotation to $\Sig_{-1}$ we
obtain a band with the same support: but this time sliding the
band across the orientation disk we obtain exactly $\Sig_1$, and
this ends the proof.~$\triangle$

 \e{dimostraz}

This theorem shows that it is not possible to extend the
definition of number of half twists to generic bands in
non-orientable $M$, since they can be equivalent up to regular
homotopy plus reparametrization or simply up to regular homotopy.

\section{Computation of cobordism.}

 We
adapt this proof to the non-orientable case.

\subsection{Statement and scheme of proof.}

Consider an immersion $f$ of a compact closed surface $F$, non
necessarily connected nor orientable, in a non-orientable
connected 3-manifold non necessarily compact nor closed. Up to
regular homotopy $f$ can be considered generic: in this case this
means it has a finite number of curves of double points and a
finite number of triple points. We define 3 invariants associated
to $f$: let $H_f$ be the homology class represented by $f$ in
$H_2(M,\zd)$; let $\del_f$ be the homology class represented by
the locus of double points of $f$, that is, the points of $M$
having 2 or 3 preimages, in $H_1(M,\zd)$; let $n_f$ be the Euler
characteristic of $F$ modulo 2 (if $F$ is not connected we
consider the sum of the Euler characteristics of its connected
components).

\be{lem}These functions are invariant up to cobordism.

\e{lem}

\be{dimostraz}Homology between cobordant immersions is provided by
the immersion of the cobordism, so that $H_{-}$ is well-defined up
to cobordism. Homology between the loci of double points of
cobordant immersions is provided by the locus of double points of
the immersion of the cobordism, so that $\del_-$ is well-defined
up to cobordism. As long as $n_-$ is concerned, recall that the
Euler class modulo 2 is the total invariant of classical cobordism
of surfaces, so that in particular it is invariant up to cobordism
of immersions.~$\triangle$

\e{dimostraz}

If $Y$ is a cobordism class of immersions we can then define
$H_Y$, $\del_Y$ and $n_Y$ as the invariants associated to any
generic representative of $Y$. We will prove that the function
\[
\be{array}{rccc} \Psi:&N_2(M)&\longrightarrow &H_2(M,\zd)\times
H_1(M,\zd)\times\zd\\ & Y&\mapsto&(H_Y,\del_Y,n_Y)\e{array}
\]
is a total invariant. First remark that
\[
\Psi(Y+Y')=(H_Y+H_{Y'},\del_Y+\del_{Y'}+H_Y\cdot H_{Y'},
n_Y+n_{Y'}),
\]
where $\cdot$ is the bilinear form of intersection; it is
straightforward to check that the composition law
\[
(H,\del,n)*(H',\del',n')=(H+H',\del+\del'+H\cdot H',n+n')
\]
makes $H_2(M,\zd)\times H_1(M,\zd)\times\zd$ a commutative group,
and $\Psi$ is clearly an homomorphism for this group structure.

\be{teo}\label{cobordismo}$\Psi$ is an isomorphism of groups
between $N_2(M)$ and $(H_2(M,\zd)\times H_1(M,\zd)\times\zd,*)$.

\e{teo}

\be{dimostraz}We are left to prove that $\Psi$ is bijective.

We first see surjectivity. Given a triple $(H,\del,n)$ consider:
an embedding $f$ that represents $H$; an embedded circle $K$
representing $\del$, its tubular neighborhood $N$ and the
immersion $h$ obtained by adding a kink along any longitude of the
inclusion of $\partial N$ (remark this is the immersion of a torus
or of a Klein bottle according to $w_1(M)(\del_Y)$ being 0 or 1);
if $n_f=n$ then $\Psi(f+h)=(H,\del,n)$; if $n_f\neq n$ then
consider $g=\phi\circ\ga$, \ga\ being a generator of $N_2(\R^3)$
and $\phi$ a diffeomorphism of $\R^3$ to a ball $B$ in $M$: then
$\Phi(f+h+g)=(H,\del,n)$; and this proves surjectivity.

Consider an immersion $f$ of  $F$ in $M$ and let $l$ be  a closed
circle in $F$; we denote by $q_f(l)$ the number of half twists of
the band given by the restriction of $f$ to a tubular neighborhood
of $l$ in $F$, possibly perturbed by a regular homotopy. Remark
that, since $M$ is non-orientable, this is only defined when $l$
is orientable in $F$ and $f(c)$ is orientable in $M$.
 To prove injectivity we show in lemma~\ref{scomposizione}, adapting the argument in
\cite{BeSSPS}, that any generic representative of a class $Y$ can
be deformed via surgeries to a representative $f+h+g$ such that:
$f$ is an embedding; $h$ is the immersion obtained from the
inclusion  of the boundary of a tubular neighborhood of a curve
representing $\del_Y$ by adding a kink along a longitude, possibly
such that $q_h(l)=0$ if $\del_Y$ is orientable; $g$ is an
immersion contained in a ball in $M$. This decomposition splits
$\Psi$ in its 3 components, that is
$\Psi(Y)=(H_f,\del_h,n_f+n_g)$, so that if $\Psi(Y)=\Psi(Y')$ we
are left to show that \be{eqnarray}
\label{embedding}H_f=H_{f'}&\Rightarrow f\sim_c f'\\
\label{ptidoppi}\del_h=\del_{h'}&\Rightarrow h\sim_c h'\\
\label{disco}n_g=n_{g'}&\Rightarrow g\sim_c g', \e{eqnarray} where
$\sim_c$ means cobordism relation. That (\ref{embedding}) is true
is the content of lemma~\ref{lembedding}; that~(\ref{ptidoppi}) is
true is the content of lemmas~\ref{lptidoppior}
and~\ref{lptidoppinonor}; that~(\ref{disco}) is true is the
content of lemma~\ref{ldisco}.~$\triangle$

\e{dimostraz}

\subsection{Decomposition lemma.}

\be{lem}\label{scomposizione} Any generic representative of a
class $Y\in N_2(M)$ can be deformed via surgeries to a
representative $f+h+g$ such that: $f$ is an embedding; $h$ is the
immersion obtained from the inclusion  of the boundary of a
tubular neighborhood of a curve representing $\del_Y$ by adding a
kink along a longitude, possibly such that $q_h(l)=0$ if $\del_Y$
is orientable; $g$ is an immersion contained in a ball in
$M$.\e{lem}

 \be{dimostraz} Let $Y$ be a
cobordism class, and let $f_1$ be a generic representative of $Y$:
it is an immersion with a finite number of curves of double points
and a finite number of triple points. We call $C(M)$ the subgroup
of $N_2(M)$ given by classes which admit a representative immersed
in a disk of $M$.

First we eliminate triple points. To do so call \ga\ the right
immersion of Boy in $\R^3$; recall \ga\ is an immersion with a
single triple point; consider the connected sum of $f_1$ and \ga\
in a chart of $M$; it is possible to deform by regular homotopy
this immersion  to an immersion with one triple point less (remark
that \ga\ in a chart of $M$ is cobordant to its inverse, see
proposition~\ref{C(M)}). By recursively repeating this
construction we obtain an immersion $f_2$ without triple points
and cobordant to $f_1$ up to an element of $C(M)$.

\begin{figure}
\begin{picture}(270,150)
\put(30,30){ \includegraphics{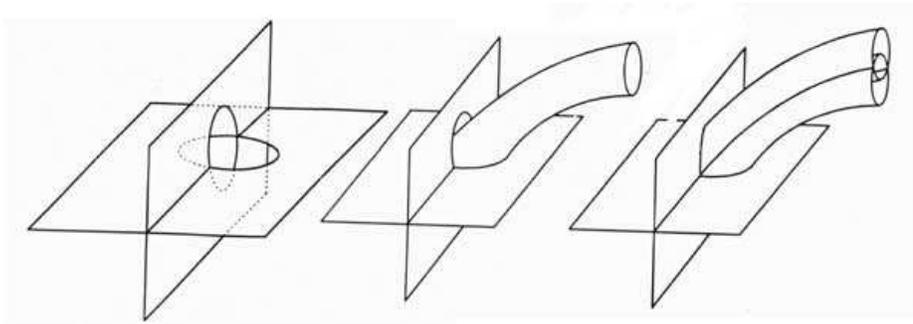} }
\end{picture}
\caption{Connecting two curves of double points.}\label{conn}
\end{figure}

The immersion $f_2$ has a finite number of curves of double
points. Take two such curves, connect them with a path, then
substitute the path with a couple of tubes as in
figure~\ref{conn}.

By recursively repeating this operation we are then left with an
immersion $f_3$ representing $Y$ up to an element of $C(M)$ and
with a single curve $K$ of double points; clearly $K$ represents
$\del_Y$.

The intersection of the image of $f_3$ with a tubular neighborhood
of $K$ is a bundle on \su\ with fiber a figure $X$, orientable or
non-orientable according to the orientability of $\del_Y$. If we
number orderly the edges of the figure $X$ we can identify the its
group of isometries with a subgroup of ${\cal S}_4$; in this
framework we can reduce the structure group of our fiber bundle to
one and only one of the following:
\begin{enumerate}
\item $G_0=1$;
\item $G_1=<(1234)>$;
\item $G_2=<(13)(24)>$;
\item $G_3=<(1432)>$;
\item  $G_4=<(12)(34)>$;
\item  $G_5=<(24)>$;
\item  $G_6=<(14)(23)>$;
\item  $G_7=<(13)>$;
\e{enumerate} call $L_0,\dots,L_7$ the corresponding bundles;
remark the first 4 are orientable, the last 4 are not.

Now remark that the groups with even index act also on the figure
8 obtained by connecting with arcs edges 1 and 4 and edges 2 and
3; this implies that in the fiber bundles with even index it is
possible to substitute the fiber, and so we obtain respectively
immersions of a torus in a solid torus, of a Klein bottle in a
solid torus, of  a torus in a solid Klein bottle and of a Klein
bottle in a solid Klein bottle.

Go back to the fiber bundle on the curve of double points. If it
is isomorphic to  $L_2$ or to $L_4$ we immerse along a curve
homotopically trivial in $M$ the fiber bundle in figure 8's
obtained from $L_2$ by substitution of fiber  and we connect its
curve of double points with $K$: the curve of double points of the
new immersion has tubolar neighborhood isomorphic to $L_0$ or to
$L_6$, respectively, and still represents $Y$ up to elements of
$C(M)$. We can then perform a Rohlin surgery as in
figure~\ref{rohlin},
\begin{figure}
\begin{picture}(270,150)
\put(30,30){ \includegraphics{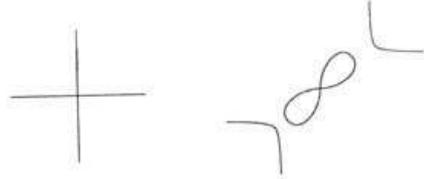} }
\end{picture}
\caption{Rohlin surgery.}\label{rohlin}
\end{figure}
and obtain from one side an embedding $f$ and from the other a map
$h$ obtained from the inclusion of the boundary of a tubular
neighborhood of a curve representing $\del_Y$ by adding a kink
along a longitude $l$; and the sum of the class of $f$ and $h$
differs from $Y$ by an element $g\in C(M)$. So if $\del_Y$ is
non-orientable we are done. If $\del_Y$ is orientable instead we
have to consider the case that $q_h(l)=2$. If it is so consider in
a disk of $M$ an immersion $h'$ obtained from the standard  by
adding a kink along a  longitude $l'$ such that $q_{h'}(l')=2$; if
we connect the curve of double points of our immersion to the
curve of double points of $h'$ in the usual manner we obtain an
immersion we call again $h$ that satisfies $q_h(l)=0$.

Now go to the case when the tubular neighborhood of $K$ is
isomorphic to a $L_i$ with $i$ odd. We have to consider another
auxiliary immersion: take on $[-1,1]$ the fibration with fibers as
in figure~\ref{fibre},
\begin{figure}
\begin{picture}(270,150)
\put(30,30){ \includegraphics{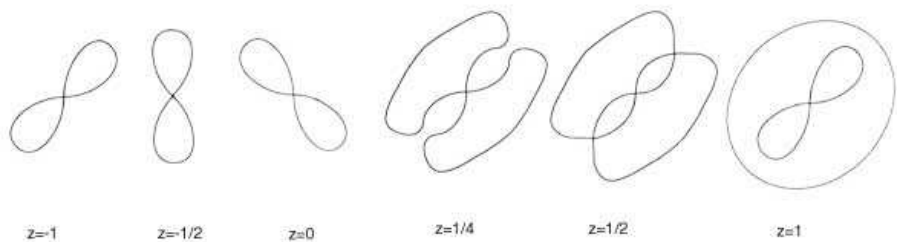} }
\end{picture}
\caption{Fibers on $[-1,1]$.}\label{fibre}
\end{figure}
and complete it to the immersion of a closed surface in a disk of
$M$ by identifying the two figure 8's at the two edges (without
torsions) and close the hole of the fiber on 1 with a 2-disk $D$;
the resulting immersion has two curves of double points meeting in
a triple point in $D$; the first has tubular neighborhood
isomorphic to $L_1$; by connecting this curve to $K$ in a new
curve $K'$ we get an immersion $f_4$, still representing $Y$ up to
elements of $C(M)$, with the tubular neighborhood of $K'$
isomorphic to an $L_i$ with $i$ even; and repeating the previous
construction we can assume $i$ equal to 0. The other curve of
double points of $f_4$, say $K''$, is contained in  $D$ and is
shown in figure~\ref{otto}.
\begin{figure}
\begin{picture}(270,150)
\put(30,30){ \includegraphics{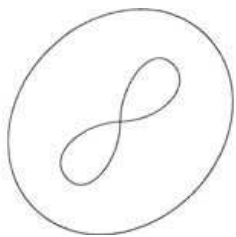} }
\end{picture}
\caption{The curve of double points $K''$ of the immersion $f_4$.
.}\label{otto}
\end{figure}
Perform Rohlin surgery on $K'$: on $D$ the result is shown in
figure~\ref{chir}:
\begin{figure}
\begin{picture}(270,150)
\put(30,30){ \includegraphics{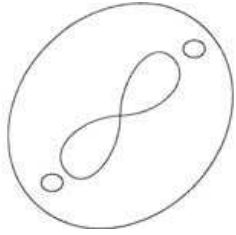} }
\end{picture}
\caption{Result of Rohlin surgery on the disk $D$.}\label{chir}
\end{figure}
to the two curves so created one can again perform Rohlin surgery
and obtain elements of $C(M)$. So, up to again applying the
previous construction when $\del_Y$ is orientable, we have a
decomposition satisfying the required properties.~$\triangle$

\e{dimostraz}

\subsection{Other lemmas.}
\be{lem}\label{lembedding} Let $f$ and $f'$ be two embeddings of
compact closed surfaces $F$ and $F'$, respectively, in the same
non-orientable 3-manifold $M$. Then $f\sim_c f'$ if and only if
$H_f=H_{f'}$.

\e{lem}

\be{dimostraz}The ``if'' part is invariance of $H_-$ up to
cobordism, and has already be proved.

For the ``only if'' part let $H$ be $H_f=H_{f'}$. If $M$ is
compact and without boundary $f$ and $f'$ are essentially the loci
of zeros of two transverse sections $s$ and $s'$ of the same line
bundle $L$, Poincar\'e dual to $H$, and this allows to construct
the cobordism this way: consider a smooth homotopy $s_t$ between
$s$ and $s'$, that be constant for $t$ close to 0 and to 1;
consider on $M\times I$ the line bundle pull-back of $L$, and
consider $s_t$ as a section of this line bundle: $s_t$ is then
transverse and its locus of zeros is the desired cobordism. If $M$
is not compact or closed in order to apply Poincar\'e duality
consider a compact 3-manifold $N$ with boundary, contained in $M$
and containing in its interior a 3-chain that bounds $f+f'$;
consider the compact closed 3-manifold $\bar{N}$ obtained by
gluing to $N$ a second copy of $N$ itself; again $f$ and $f'$
represent the same element $H$ in the homology of $\bar{N}$, so
they are loci of zeros of transverse sections $s$ and $s'$ of the
same line bundle $L$, Poincar\'e dual (in $\bar{N}$) to $H$; we
can assume that $s$ and $s'$ coincide outside a compact contained
in $N$; we construct as before a smooth homotopy $s_t$, that we
pretend relative to a compact containing the complementary of $N$;
then the construction runs as before, and gives a cobordism which
is in fact contained in $N\times I$, hence in $M\times
I$.~$\triangle$

\e{dimostraz}


We now prove (\ref{ptidoppi}). We first introduce some notation.
We denote \cl{f}\ the class of $f$ up to reparametrizations of
$F$, that is, $g$ belongs to \cl{f}\ if there exists a
diffeomorphism $\phi$ of $F$ such that $g=f\circ\phi$; remark that
regular homotopy equivalence relation is well defined up to
reparametrization, and that immersions that differ by a
reparametrization are cobordant (see \cite{BeSSPS}, pages 657--8).
Finally denote by $\cl{f}*K$ the reparametrization class of
$f*f^{-1}(K)$, that is well defined. The key of the proof are
propositions 10.5 and 10.7 of \cite{BeSSPS}, that we  adapt to a
non-orientable situation:

\be{lem}Let $f$ be an immersion of $F$ in $M$. Let $K$ be a closed
circle in $f(F)$ such that $q_{f}(f^{-1}(K))=0$. Then:
 \be{enumerate}
\item if $K$ is trivial in $\pi_1(M)$ then \cl{f}\ and
$\cl{f}*K$ are regularly homotopic;
\item if $K$ is trivial in  $H_1(M,\zd)$ then \cl{f}\ and
$\cl{f}*K$ are cobordant. \e{enumerate} \e{lem}

\be{dimostraz}1. Let $c=f^{-1}(K)$; $c$ is orientable in $F$. Let
$\phi$ be a twist of Dehn of $F$ along $c$, we prove that $f*c$
and $f\circ \phi$ are regularly homotopic. We make direct use of
the result of \cite{HaHIST}.

First remark that, since $K$ is homotopically trivial,
$f\circ\phi$ is homotopic to $f$, that on its side is homotopic to
$f*c$. We can then consider on $f*c$ and $f\circ\phi$ the
correspondence $C_f$, as defined in the proof of
theorem~\ref{teoisotropia}; we see that
$C_f(f*c)=C_f(f\circ\phi)$.

Consider an embedding of the standard torus  in $\R^3$ with the
\spin\ structure induced by restriction of the unique \spin\
structure on $\R^3$   as a tubular neighborhood of $K$, in such a
way that the \spin\ structure of either of the two preimages of
$K$ in \tM\ is respected. The ipothesis $q_f(c)=0$ implies that
this embedding can be chosen in such a way that $f$ restricted to
a tubular neighborhood of $c$ can be read as the embedding in the
standard torus  of a band with no half twists. In the same model
then $f\circ\phi$ and $f*c$ are read as shown in
figure~\ref{cilindri}, and far from $K$ they coincide with $f$.
\begin{figure}
\begin{picture}(270,150)
\put(30,30){ \includegraphics{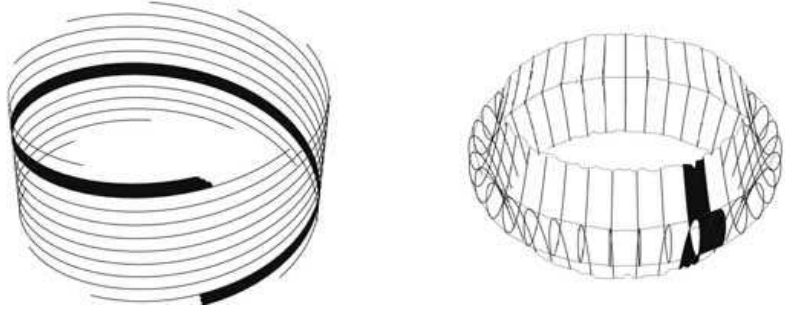} }
\end{picture}
\caption{Local model of $f\circ\phi$ and of
$f*c$.}\label{cilindri}
\end{figure}

Now take a representative $d$ of an element $H_1(F,\zd)$,
transverse to $c$: following again the proof of the main theorem
of \cite{HaHIST} one shows that if $d\cdot c=0$ then
$C_f(f\circ\phi)(d)=C_f(f*c)(d)=0$, and if $d\cdot c=1$ then
$C_f(f\circ\phi)(d)=C_f(f*c)(d)=1$; and this ends the proof.

2. The proof of \cite{BeSSPS} can be used, by means
of~1.$\triangle$
 \e{dimostraz}


We remark that, since immersions that differ by a
reparametrization are cobordant, the second part of the lemma
implies:

\be{coroll}\label{condizione}If $K$ is in the image of an
immersion $f$ and is trivial in $H_1(M,\zd)$ then $f$ is cobordant
to $f*f^{-1}(K)$.~$\triangle$ \e{coroll}

 We are now ready for the
lemmas that prove (\ref{ptidoppi}).

\be{lem}\label{lptidoppior}Let $\del\in H_1(M,\zd)$ such that
$w_1(M)(\del)=0$; let $K$ and $K'$ be closed circles both
representing \del; let $i$ and $i'$ be the inclusions of the
boundaries of tubular neighborhoods of $K$ and $K'$, respectively.
Let $h$ and $h'$ be two immersions obtained from $i$ and $i'$ by
adding kinks along longitudes $l$ and $l'$ such that
$q_i(l)=q_{i'}(l')=0$. Then $h\sim_c h'$.

\e{lem}

\be{dimostraz}This is proven in \cite{BeSSPS}, page 672. We review
their argument. Remark that $i$ and $i'$ are immersions of tori,
representing the identity of the semi-group structure; their
connected sum $\bar{i}=i\sharp i'$ still represents the identity.
The sum of the cobordism classes of $h$ and $h'$ is represented by
their connected sum $\bar{h}=h\sharp h'$; remark that $\bar{h}$ is
obtained from $\bar{i}$ by adding a kink along a curve $\bar{l}$
homologous to $l+l'$, so that
$q_{\bar{i}}(\bar{l})=q_i(l)+q_{i'}(l')=0$. Moreover the homology
class of $\bar{i}(\bar{l})$ is $2\del$, that is 0, in
$H_1(M,\zd)$. So we can apply corollary~\ref{condizione} and say
that $\bar{h}$ is cobordant to $\bar{i}$, that is the identity; so
$h'$ belongs to the opposite of the class of $h$.

But now do the same construction with two copies of $h$, and
conclude that $h$ itself belongs to its own opposite class, hence
the claim.~$\triangle$

 \e{dimostraz}

To settle the non-orientable case we first need a remark:

\be{lem}\label{kleininklein} Let $K$ be a curve non-orientable in
$M$, let $i$ be the inclusion of the boundary of a tubular
neighborhood $N$ of $K$, let $l$  and $m$ be respectively a
longitude and the meridian of $\partial N$; then
\[
i*l\sim_c i*(l+m).
\]
\e{lem}

\be{dimostraz}The homotopy class of $i$ is odd: in fact $i$ can be
homotopically deformed to the inclusion of $K$, this can slide
across a surface representing $w_1(M)$ and then going back to $i$,
and these three steps give a odd self-homotopy of $i$. But now $m$
represents the dual to $w_1(\partial N)$, so that by
proposition~\ref{teoisotropia} adding a kink along $m$ doesn't
change regular homotopy class, hence cobordism class.~$\triangle$
 \e{dimostraz}

\be{lem}\label{lptidoppinonor}Let $\del\in H_1(M,\zd)$ such that
$w_1(M)(\del)\neq 0$; let $K$ and $K'$ be closed circles both
representing \del; let $i$ and $i'$ be the inclusions of the
boundaries of tubular neighborhoods of $K$ and $K'$, respectively.
Let $h$ and $h'$ be two immersions obtained from $i$ and $i'$ by
adding kinks along longitudes $l$ and $l'$. Then $h\sim_c h'$.

\e{lem}

\be{dimostraz}The inclusions $i$ and $i'$ are immersions of Klein
bottles, both representing the identity of the cobordism
semi-group. Their connected sum $\bar{i}=i\sharp i'$ still
represents the identity. The sum of the cobordism classes of $h$
and $h'$ is represented by their connected sum $h\sharp h'$, which
is an immersion regularly homotopic to the one obtained from
$\bar{i}$ by adding a kink along a curve $\bar{l}$ homologous to
$l+l'$; this is cobordant to the immersion obtained from $\bar{i}$
by adding a kink along the curve $\bar{l}'$ homologous to
$l+l'+m$, $m$ being the meridian of the tubular neighborhood of
$K$, because of lemma~\ref{kleininklein}. But both of $\bar{l}$
and $\bar{l}'$ represent $2\del$, hence 0, in $H_1(M,\zd)$, and
one of $q_{\bar{i}}(\bar{l})$ or $q_{\bar{i}}(\bar{l}')$ must be
0, see figure~\ref{bibiklein}.
\begin{figure}
\begin{picture}(270,150)
\put(30,30){ \includegraphics{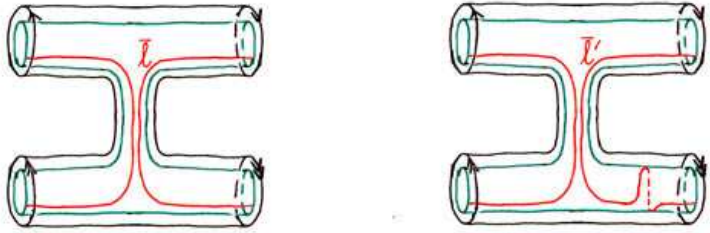} }
\end{picture}
\caption{The curves $\bar{l}$ and $\bar{l}'$ on the image of the
embedding $\bar{i}$. }\label{bibiklein}
\end{figure}
So we can apply corollary~\ref{condizione}, and obtain that
$\bar{h}$ is cobordant to $\bar{i}$, that is the identity class.
This implies that $h'$ belongs to the opposite class to $h$.

But now applying the same construction to two copies of $h$ we
obtain that $h$  itself belongs to its own opposite class, hence
the claim.~$\triangle$

\e{dimostraz}

We are left now to the part that is in some way more
characteristic of the non-orientable case, that is, the point that
gives account of the fact that the cobordism groups of
non-orientable 3-manifolds are ``smaller" than the groups of the
orientable case. The key is again the fact that in a
non-orientable environment there is a kind of isotopy that doesn't
exist in an orientable manifold, that is, an isotopy that reverses
orientation.

\be{prop}\label{C(M)}Let $C(M)$ be the subgroup of $N_2(M)$ given
by classes which admit a representative immersed in a disk of $M$.
$C(M)$ is isomorphic to \zd.

\e{prop}

\be{dimostraz}Consider a coordinate chart $(U,\phi)$ of $M$. The
diffeomorphism $\phi$ induces a homomorphism
\[
\be{array}{cccl}\phi^*:& N_2(\R^3)&\longrightarrow&C(M)\\
&f&\mapsto& \phi\circ f\e{array}
\]
that is evidently surjective. Recall that $N_2(\R^3)$ is a cyclic
group of order 8 generated by the right immersion of Boy  of
$\Pro^2$; we call \ga\ this immersion. This implies that $C(M)$ is
a cyclic group generated by $\phi^*(\ga)$, that can't be trivial
since $\Pro^2$ generates the cobordism group of surfaces.

Let now $f_t$ be a self-isotopy of $B$ that reverses the
orientation, for example that crosses once a surface representing
$w_1(M)$. Recall that the inverse of $\ga$ in $N_2(\R^3)$ is given
by composition with a reflection in a plane, so that
$f_1\circ\phi^*(\ga)$ is the image via $\phi^*$ of the canonical
representative of the opposite class to \ga, hence represents the
opposite class to $\phi^*(\ga)$. But $f_1\circ\phi^*(\ga)$, being
isotopic to $\phi^*(\ga)$, is cobordant to it, so that \ga\
belongs to its own opposite class, hence has order 2.~$\triangle$

\e{dimostraz}

\be{lem}\label{ldisco}Let $g$ and $g'$ be immersions whose image
is contained in a disk of $M$; then $g\sim_c g'$ if and only if
$n_g=n_{g'}$.

\e{lem}

\be{dimostraz}The ``if'' part is invariance of $n_-$ up to
cobordism, and has already be proved. For the ``only if" it is
enough to remark that $n_-$ realizes the homomorphism of
proposition~\ref{C(M)}.~$\triangle$

\e{dimostraz}

\vspace{1.5cm}

\nnd Rosa Gini, {\tt gini@dm.unipi.it},\\ Dipartimento di Matematica
``Leonida Tonelli'',\\ via Filippo Buonarroti 2, \\ I--56127 Pisa,
Italy.


\begin{thebibliography}{Wel66}

\bibitem[Bro72]{BroGKI}
E.H. Brown.
\newblock Generalizations of the {K}ervaire invariant.
\newblock {\em Ann. Math.}, 95:368--383, 1972.

\bibitem[BS95]{BeSSPS}
R.~Benedetti and R.~Silhol.
\newblock ${S}pin$ and ${P}in^-$ structures, immersed and embedded surfaces and
  a result of {S}egre on real cubic surfaces.
\newblock {\em Topology}, 34:651--678, 1995.

\bibitem[GM86]{GuMETR}
L.~Guillou and A.~Marin.
\newblock Une estension d'un th\'eor\`eme de {R}ohlin sur la segnature.
\newblock In L.~Guillou and A.~Marin, editors, {\em \`A la recherche de la
  topologie perdue}, pages 97--117. Birkh\"auser, Basel, 1986.

\bibitem[HH85]{HaHIST}
J.~Hass and J.~Hughes.
\newblock Immersions of surfaces in 3-manifolds.
\newblock {\em Topology}, 24:97--112, 1985.

\bibitem[Hir59]{HirIM}
M.~W. Hirsch.
\newblock Immersions of manifolds.
\newblock {\em Trans. Amer. Math. Soc.}, 93:242--277, 1959.

\bibitem[KT89]{KiTPSL}
R.~C. Kirby and L.~R. Taylor.
\newblock ${P}in$ structures on low-dimensional manifolds.
\newblock In S.~K. Donaldson and C.~D. Thomas, editors, {\em Geometry of
  low-dimensional manifolds}, pages 177--241. Cambridge University Press,
  Cambridge, 1989.

\bibitem[Mil62]{MilSSM}
J.~Milnor.
\newblock ${S}pin$ structures on manifolds.
\newblock {\em L'enseignement Ma\-th\'e\-ma\-ti\-que}, 8:198--203, 1962.

\bibitem[Pin85]{PinRHC}
U.~Pinkall.
\newblock Regular homotopy classes of immersed surfaces.
\newblock {\em To\-po\-lo\-gy}, 24:421--434, 1985.

\bibitem[Wel66]{WelCGI}
R.~Wells.
\newblock Cobordism groups of immersions.
\newblock {\em To\-po\-lo\-gy}, 5:281--294, 1966.

\end{thebibliography}
\end{document}